\documentclass[11pt]{amsart}
\usepackage{amsmath}
\usepackage{amssymb,array,stmaryrd}
\usepackage{mathtools}
\usepackage{enumerate}
\usepackage[all,cmtip]{xy}

\begin{document}

\newtheorem{thm}{Theorem}[section]
\newtheorem*{thm*}{Theorem}
\newtheorem{lem}[thm]{Lemma}
\newtheorem*{lem*}{Lemma}
\newtheorem{prop}[thm]{Proposition}
\newtheorem*{prop*}{Proposition}
\newtheorem{cor}[thm]{Corollary}
\newtheorem{defn}[thm]{Definition}
\newtheorem*{conj*}{Conjecture}
\newtheorem{conj}[thm]{Conjecture}
\theoremstyle{remark}
\newtheorem*{remark}{Remark}
\newtheorem*{question}{Question}

\numberwithin{equation}{section}

\newcommand{\Z}{{\mathbb Z}} 
\newcommand{\Q}{{\mathbb Q}}
\newcommand{\R}{{\mathbb R}}
\newcommand{\C}{{\mathbb C}}
\newcommand{\N}{{\mathbb N}}
\newcommand{\NC}{{\mathcal N}}
\newcommand{\PR}{{\mathcal P}}
\newcommand{\FF}{{\mathbb F}}
\newcommand{\fq}{\mathbb{F}_q}
\newcommand{\fqtm}{\mathbb{F}_q[t]^{<m}}
\newcommand{\feq}{\overline{\mathbb F}_q}
\newcommand{\slo} {\operatorname{sl}}
\newcommand{\disc}{\operatorname{disc}}
\newcommand{\tr}{\operatorname{tr}}
\newcommand{\sktr}{\widetilde{\operatorname{tr}}}
\newcommand{\tPsi}{\widetilde\Psi}
\newcommand{\tPhi}{\widetilde\Phi}
\newcommand{\No}{\operatorname{N}}
\newcommand{\mingap}{\delta_{\min}^{(\alpha)}(N)}
\newcommand{\eps}{\varepsilon}
\newcommand{\OK}{\mathcal{O}_K}
\newcommand{\cnj}{\overline}

\DeclarePairedDelimiter{\floor}{\lfloor}{\rfloor}

 \title[Evenly divisible rat. approx. of quad. irrationalities]
 {Evenly divisible rational approximations of quadratic irrationalities}
  \author{Dan Carmon}
 \address{Raymond and Beverly Sackler School of Mathematical Sciences,
Tel Aviv University, Tel Aviv 69978, Israel}
 \date{\today}
 \thanks{The research leading to these results has received funding from the European
Research Council under the European Union's Seventh Framework Programme
(FP7/2007-2013) / ERC grant agreement n$^{\text{o}}$ 320755.}

\begin{abstract}
 In a recent paper of Blomer, Bourgain, Radziwi\l\l\ and Rudnick \cite{BBRR},
 the authors proved the existence of small gaps between eigenvalues of 
 the Laplacian in a rectangular billiard with sides $\pi$ and $\pi/\sqrt\alpha$,
 i.e.\@ numbers of the form $\alpha m^2+ n^2$, whenever $\alpha$ is a quadratic
 irrationality of certain types. In this note, we extend
 their results to all positive quadratic irrationalities $\alpha$.  
\end{abstract}

 \maketitle
 
 \section{Introduction}
 
 In \cite{BBRR}, the authors investigated the minimal gaps between energy 
 levels of the eigenvalues of the Laplacian of a billiard in a rectangle
 with width $\pi/\sqrt\alpha$ and height $\pi$. These eigenvalues are of the
 form $\alpha m^2 + n^2$, where $m,n \ge 1$ are integers. For positive irrational
 $\alpha$, these levels belong to a simple
 spectrum $0 < \lambda_1 < \lambda_2 < \cdots $, with growth
 $\lambda_N \sim \frac{4\sqrt{\alpha}}{\pi} N$. The minimal gap function
 is then defined as 
 $$\mingap = \min(\{\lambda_{i+1}-\lambda_i : 1 \le i < N\}).$$
 
 The authors prove several lower and upper bounds on the growth rate of 
 $\mingap$, for various families and sets of irrational $\alpha$. Amongst these
 is \cite[Theorem 6.1]{BBRR}:
 \begin{thm}\label{thm 61}
  For all positive real quadratic irrationalities of the form
  \begin{equation}\label{quad form}
   \alpha = \alpha(x; a,b,\epsilon,r) =  r\cdot\left(
   \frac{x+\sqrt{x^2+4\epsilon}}{2}\right)^a
   \cdot\left(\sqrt{x^2+4\epsilon}\right)^b
  \end{equation}
  with 
  $$a \in \Z,\quad b = 0,1,\quad  x \in \Z\setminus\{0\},
  \quad \epsilon = \pm 1, \quad r \in \Q^{\times},$$
  we have $\mingap \ll_{\alpha,\eps} N^{-1+\eps}$ infinitely often, for any $\eps > 0$.
 \end{thm}
 Furthermore, the authors remark that when $b=0$ and $a$ is even, 
 they have in fact obtained the stronger result that $\mingap \ll_\alpha N^{-1}$
 for all $N$. We refer the reader to \cite{BBRR} for a more 
 detailed introduction to the problem, as well as the motivation
 for these upper bounds.
 
 Our goal in this paper is to generalize the result above to any positive quadratic 
 irrationality. Here on, we shall denote by $K = \Q[\sqrt{D}]$ -- a real 
 quadratic number field, $\OK$ -- its ring of integers, 
 and for any $\omega \in K$, we shall denote its conjugate in $K$
 by $\cnj\omega$, and its norm and trace (over $\Q$) by 
 $\No(\omega) = \omega\cnj\omega$, $\tr(\omega) =\omega + \cnj\omega$.
 In section \ref{main sect}, we prove:
 
 \begin{thm}\label{main thm}
  Let $K = \Q[\sqrt{D}]$ be any real quadratic number field, where 
  $D \in \N$ is square-free.
  Let $0 < \alpha \in K$ be a positive element of $K$. Then
  $\mingap \ll_{\alpha,\eps} N^{-1+\eps}$ infinitely often, for any $\eps > 0$.
  Furthermore, if \linebreak
  $\alpha \in \Q_{>0}\cdot(K^\times)^2$,\footnote{Given any $\alpha \in K$,
  this condition can 
  be easily checked, as it is equivalent to $\alpha$ satisfying
  $\No(\alpha) \in (\Q^\times)^2$ and $\tr(\alpha) > 0$.
  We do not use this equivalence in this paper, and leave its proof
  as a simple exercise to the reader.}
  then $\mingap \ll N^{-1}$ for all $N$. 
 \end{thm}
 
 One of the key lemmas in \cite{BBRR} relates small gaps between
 the $\lambda_i$ to finding good rational approximations to $\alpha$
 with both numerator and denominator being {\it evenly divisible} (or
 {\it strongly even divisible}): We call a sequence $m_n$ of numbers
 evenly divisible if there exist divisors $d_n \mid m_n$ such that
 $\min(d_n,m_n/d_n) \gg_\eps m_n^{1/2-\eps}$, for all $\eps > 0$. We call
 the sequence strongly evenly divisible if $\min(d_m,m_n/d_n) \gg m_n^{1/2}$.
 Then by \cite[Lemma 3.1]{BBRR}, we have:
 
 \begin{lem}\label{lem 31}
  If $\alpha > 0$ has infinitely many rational approximations
  $p_n/q_n$ satisfying 
  \begin{equation}\label{good approx}
   |q_n\alpha - p_n| \ll \frac{1}{q_n}
  \end{equation}
  with the sequences $\{p_n\}, \{q_n\}$ being evenly divisible (resp. 
  strongly evenly divisible), then $\mingap \ll N^{-1+\eps}$
  for all $\eps > 0$ (resp. $\mingap \ll N^{-1}$) infinitely often.
  
  If in addition $q_n \gg q_{n+1}$ for all $n$, then these inequalities 
  hold for all $N$. 
 \end{lem}
 
 The authors of \cite{BBRR} were able to find evenly divisible approximations
 satisfying \eqref{good approx} for quadratic irrationalities of the form
 \eqref{quad form}. However, the condition \eqref{good approx} can be 
 significantly weakened; namely, if we replace it with the condition
 
 \begin{equation}
  |q_n\alpha - p_n| \ll \frac{1}{q_n^{1-\eps}}, \text{ for all } \eps > 0
 \end{equation}
 
 The conclusion of Lemma \ref{lem 31} that $\mingap \ll N^{-1+\eps}$ for all
 $\eps > 0$ remains valid. This is essentially the content of 
 \cite[Lemma 3.2]{BBRR}, and the form of the lemma that will be more useful
 to us.
 
 \section{Divisibility properties}
 
 In this section we establish some divisibility properties
 for traces and twisted traces of powers of a quadratic algebraic integer. These
 play similar roles to the Chebyshev polynomials featured in \cite{BBRR},
 and in fact generalize them.
 
 \subsection{Divisibility properties of traces}\label{trace sec}
 Let $\omega \in \OK$ be an algebraic integer with $\tr\omega \neq 0$. Note that
 for any odd positive integer $\ell$, $\omega^\ell \in \OK$, hence $\tr(\omega^\ell)$
 is an integer. Furthermore, it is divisible by $\tr(\omega)$;
 this follows from the fact that $\frac{x^\ell + y^\ell}{x+y}$ is a 
 symmetric polynomial in $\Z[x,y]$,
 and substituting $\omega, \cnj\omega$ for $x,y$. More generally,
 for any positive integer $L$ which is odd and square-free, 
 we have the decomposition
 \begin{equation}
  \tr(\omega^L) = \Phi_L(\omega)\Psi_L(\omega),
 \end{equation}
 where $\Phi_L(\omega),\Psi_L(\omega) \in \Z$ are both integers, defined by
 \begin{equation}
  \Phi_L(\omega) = \prod_{\ell \mid L} \tr(\omega^{L/\ell})^{\mu(\ell)}, \quad
  \Psi_L(\omega) = \prod_{\ell \mid L; \ell > 1} \tr(\omega^{L/\ell})^{-\mu(\ell)}.
 \end{equation}
 Indeed, both $\Phi_L, \Psi_L$ are symmetric polynomials in $\Z[\omega,\cnj\omega]$\footnote{
 Specifically, $\Phi_L$ is the homogenized reflected $L$-th cyclotomic polynomial, and $\Psi_L$ is
 the corresponding cofactor.}, hence their values must be in $\Z$ for
 $\omega \in \OK$.
 
 In order to estimate the sizes of the factors $\Phi_L,\Psi_L$, assume further
 that $|\cnj\omega/\omega| \le 1/2$. Note that 
 $$e^{-2|x|} \le 1-|x| \le |1+x| \le 1+|x| \le e^{2|x|}$$
 for all $|x| \le 1/2$, and therefore for any $m \ge 1$, 
 $$|\tr(\omega^m)| = |\omega^m + \cnj\omega^m| = |\omega^m|\left|1+
 \left(\frac{\cnj\omega}{\omega}\right)^m\right| \in 
 |\omega^m|\cdot[e^{-2^{1-m}},e^{2^{1-m}}],$$
 hence,
 \begin{equation}
  |\tr(\omega^L)| \asymp |\omega|^L, \quad
  |\Phi_L(\omega)| \asymp |\omega|^{\phi(L)}, \quad
  |\Psi_L(\omega)| \asymp |\omega|^{L-\phi(L)},
 \end{equation}
 where the implied constants are all between $e^{-2}$ and $e^{2}$. It immediately follows that
 \begin{equation}
  |\Phi_L(\omega)| \asymp |\tr(\omega^L)|^{\tfrac{\phi(L)}{L}}, \quad
  |\Psi_L(\omega)| \asymp |\tr(\omega^L)|^{1-\tfrac{\phi(L)}{L}}.
 \end{equation}
 
 \subsection{Divisibility properties of twisted traces}
 For $\omega \in K$, we define the {\it twisted trace} of $\omega$ as
 \begin{equation}
 \sktr(\omega) := \tr(\sqrt{D}\cdot\omega) = (\omega - \cnj\omega)\sqrt{D} \in \Q.
 \end{equation}
 Suppose $\omega \in \OK$, and $\sktr(\omega) \neq 0$ (i.e. $\omega \notin \Q$).
 Then for all positive integers $\ell$,  $\sktr(\omega^\ell)$ is an integer
 and divisible by $\tr(\omega)$, since $\frac{x^\ell - y^\ell}{x-y}$ is a
 symmetric polynomial in $\Z[x,y]$ and $\sqrt{D}\omega^\ell \in \OK$.
 As above, this generalizes to any square-free $L \in \N$, by
  \begin{equation}
  \sktr(\omega^L) = \tPhi_L(\omega)\tPsi_L(\omega),
 \end{equation}
 where 
 \begin{equation}
  \tPhi_L(\omega) = \prod_{\ell \mid L} \sktr(\omega^{L/\ell})^{\mu(\ell)}, \quad
  \tPsi_L(\omega) = \prod_{\ell \mid L; \ell > 1} \sktr(\omega^{L/\ell})^{-\mu(\ell)}.
 \end{equation}
 are again integers, and whenever $|\cnj\omega/\omega| < 1/2$ and $L > 1$, we have
  \begin{equation}
  |\sktr(\omega^L)| \asymp \sqrt{D}|\omega|^L, \quad
  |\tPhi_L(\omega)| \asymp |\omega|^{\phi(L)}, \quad
  |\tPsi_L(\omega)| \asymp \sqrt{D}|\omega|^{L-\phi(L)},
 \end{equation}
 and consequentially 
 \begin{equation}
  |\tPhi_L(\omega)| \asymp {D}^{-\tfrac{\phi(L)}{2L}} |\sktr(\omega^L)|^{\tfrac{\phi(L)}{L}}, \quad
  |\tPsi_L(\omega)| \asymp {D}^{\tfrac{\phi(L)}{2L}} |\sktr(\omega^L)|^{1-\tfrac{\phi(L)}{L}}.
 \end{equation}
 
 The major difference between the divisibility properties for $\tr(\omega^L)$
 and $\sktr(\omega^L)$ is that in the second case $L$ may be even --
 and specifically, we may choose $L = 2$, which is the only number with
 $\phi(L)/L = 1/2$, capable of generating strongly evenly divisible
 sequences, rather than just evenly divisible.

 \section{Proof of Theorem \ref{main thm}}\label{main sect}
 
 We first note that we may assume $\alpha \in \OK$. Indeed, for any $\alpha \in K$
 there is some denominator $A \in \N$ such that $A\cdot\alpha \in \OK$. If we are
 then able to find (strongly) evenly divisible $p_n, q_n$ such that $A\cdot\alpha$ is
 well-approximated by $p_n/q_n$, then the sequence $A\cdot q_n$ is also (strongly)
 evenly divisible and $\alpha$ is well-approximated by $p_n/(A\cdot q_n)$.
 
 Let $\zeta \in \OK$ be an algebraic integer with $\No(\zeta) = 1$ and
 $|\zeta| > 1 > |\cnj\zeta|$. Such $\zeta$ may constructed e.g.\@ as $\zeta = x + y\sqrt{D}$,
 where $(x,y) \in \N^2$ is a non-trivial integral solution to the Pell equation
 $x^2 - Dy^2 = 1$.\footnote{We may choose $\zeta$  with $\No(\zeta) = -1$ 
 instead, in fields where such integers exist. Note that the number
 $\frac{x+\sqrt{x^2+4\epsilon}}{2}$ appearing in \eqref{quad form} 
 is always an example of an appropriate $\zeta$ for the appropriate field.}
 Let $\eps > 0$ be arbitrarily small.

 \subsection{A symmetric construction}
 Let $\{\ell_i\}_{i=1}^{\infty}, \{\ell'_i\}_{i=1}^{\infty}$ be two disjoint subsequences of the odd primes, with 
 $$\prod_{i=1}^{\infty} \left(1-\frac{1}{\ell_i}\right) = 
 \prod_{i=1}^{\infty} \left(1-\frac{1}{\ell'_i}\right) = \frac12,$$
 the existence of which easily follows from $\prod_p (1-\tfrac{1}{p}) = 0$. 
 For a fixed $t \in \N$ write 
 $$L = \prod_{i=1}^t \ell_i,\quad L' = \prod_{i=1}^t \ell'_i$$
 and suppose $t = t(\eps)$ is sufficiently large so that
 $\frac{\phi(L)}{L},\frac{\phi(L')}{L'} \in (\frac12,\frac12 + \eps)$.
 
 Let $M = M(\eps) \in \N$ be the smallest positive integer such that $M+1$ is divisible by $L$ and 
 $M$ is divisible by $L'$, and write $M+1 = m_1 L$, $M = m_2 L'$. 
 Define $N = nLL'$, where $n = n(\eps)$ is any sufficiently large number 
 such that the following inequalities are satisfied:
 \begin{align}
  \label{ineq for tr}
  |\zeta/\cnj\zeta|^{nL'} & > 2|\cnj\alpha/\alpha|^{m_1}, \quad
  |\zeta/\cnj\zeta|^{nL} > 2|\cnj\alpha/\alpha|^{m_2}, \\
  \label{ineq for C}
  |\zeta|^{\eps N} & > 4|\alpha - \cnj\alpha||\No(\alpha)|^M|\alpha|^{-\eps M}.
 \end{align}
 Such $n$ clearly exist, as $\alpha$ is fixed, $m_1$, $m_2$, $L$ and $L'$ depend only on $\eps$,
 and $|\zeta/\cnj\zeta| > |\zeta| > 1$, by the definition of $\zeta$. Finally, we define
 \begin{align}
  P &= P(\eps) = \tr(\alpha^{M+1}\zeta^{N}) = \tr((\alpha^{m_1}\zeta^{nL'})^{L}), \\
  Q &= Q(\eps) = \tr(\alpha^{M}\zeta^{N}) = \tr((\alpha^{m_2}\zeta^{nL})^{L'}).
 \end{align}
 
 We first show that $P/Q$ is a good approximation to $\alpha$. Note that
 \eqref{ineq for tr} implies $|Q|$ is approximately $|\alpha^M \zeta^N|$, or,
 more explicitly,
 \begin{equation} \label{qaz rel}
  |Q|/|\alpha^M \zeta^N| \in [1-2^{-L'},1+2^{-L'}] \subset [1/2,2].
 \end{equation}
 It then follows that
 $$|\alpha Q - P| = |\alpha-\cnj\alpha||\cnj\alpha^M\cnj\zeta^N| = 
 \frac{|\alpha-\cnj\alpha||\No(\alpha)|^M|\No(\zeta)|^N}{|\alpha^M \zeta^N|} 
 \le \frac{C(\alpha,\eps)}{|Q|} \le \frac{1}{|Q|^{1-\eps}},
 $$
 where $C(\alpha,\eps) = 2|\alpha-\cnj\alpha||\No(\alpha)|^{M(\eps)}$ 
 depends only on $\alpha$ and $\eps$. 
 The two final steps then follow immediately from \eqref{qaz rel}, \eqref{ineq for C}
 and $|\No(\zeta)| = 1$.
 
 Next, we show that $P$ and $Q$ are evenly divisible. By their definitions,
 the inequalities \eqref{ineq for tr} and
 section \ref{trace sec}, it is evident that the factorization 
 $P = \Phi_L(\alpha^{m_1}\zeta^{nL'})\Psi_L(\alpha^{m_1}\zeta^{nL'})$
 satisfies 
 $$\min(|\Phi_L(\alpha^{m_1}\zeta^{nL'})|,|\Psi_L(\alpha^{m_1}\zeta^{nL'})|)
 \asymp \min(|P|^{\tfrac{\phi(L)}{L}},|P|^{1-\tfrac{\phi(L)}{L}}) 
 \gg |P|^{1/2 - \eps},
 $$ 
 and similarly $Q = \Phi_{L'}(\alpha^{m_2}\zeta^{nL})\Psi_{L'}(\alpha^{m_2}\zeta^{nL})$
 with
 $$\min(|\Phi_{L'}(\alpha^{m_2}\zeta^{nL})|,|\Psi_{L'}(\alpha^{m_2}\zeta^{nL})|)
 \gg |Q|^{1/2 - \eps}.
 $$ 
 We conclude that for any sequence $\eps_k \ssearrow 0$, the sequences 
 $P(\eps_k),Q(\eps_k)$ are evenly divisible and provide good approximations
 to $\alpha$, which allows us to conclude via Lemma \ref{lem 31}.
 
 \subsection{Alternative constructions}
 Instead of defining $P$ and $Q$ via traces of integral powers, 
 we may use twisted traces instead, and define
 $$P = \sktr((\alpha^{m_1}\zeta^{nL'})^L),\quad
 Q = \sktr((\alpha^{m_2}\zeta^{nL})^{L'}).$$
 Repeating the same computations for the same values of $L,L'$
 would lead to the same estimates as in the first construction.
 However, since we are now using skew-traces, we may define
 $L$ (resp. $L'$) to be equal to $2$ instead of a product of odd primes,
 which will then correspond to $P$ (resp. $Q$) having two factors of size
 $\gg_{D} |P|^{1/2}$, i.e.\@ the sequence $P(\eps_k)$ (resp. $Q(\eps_k)$) will be
 strongly evenly divisible. However, we cannot in general 
 do this for both $P$ and $Q$ simultaneously, since it is impossible for both
 $M+1$ and $M$ to be even.
 
 Suppose now that $\alpha \in \Q_{>0}\cdot(K^\times)^2$. By multiplying by
 a proper natural denominator $A$, we may assume that $\alpha \in \OK^2$,
 i.e. $\alpha = \beta^2$ for some $\beta \in \OK$. Now, for all sufficiently
 large $n$, define
 \begin{align}
  P_n &:= \sktr(\alpha\zeta^{2n}) = \sktr((\beta\zeta^n)^2) = \sktr(\beta\zeta^n)\cdot\tr(\beta\zeta^n), \\ 
  Q_n &:= \sktr(\zeta^{2n}) = \sktr((\zeta^n)^2) = \sktr(\zeta^n)\cdot\tr(\zeta^n).
 \end{align}

 The sequences $P_n$ and $Q_n$ are therefore both strongly evenly divisible, since
 for large $n$ we have
 \begin{equation}
  |\sktr(\beta\zeta^n)| \asymp \sqrt{D}|\beta \zeta^n| \asymp \sqrt{D}|\tr(\beta \zeta^n)|, \\
  \end{equation}
  \begin{equation}
  |\sktr(\zeta^n)| \asymp \sqrt{D}|\zeta^n| \asymp \sqrt{D}|\tr(\zeta^n)|.
 \end{equation}
 Furthermore, their ratios approximate $\alpha$ well, as
 $$|Q_n\alpha - P_n| = |\alpha - \cnj\alpha||\cnj\zeta|^{2n} =
 \frac{|\alpha - \cnj\alpha|}{|\zeta|^{2n}} \asymp 
 \frac{|\alpha - \cnj\alpha|\sqrt{D}}{|Q_n|}, $$
 and the denominators grow geometrically -- $|Q_{n}| \asymp |\zeta^{-2}| |Q_{n+1}| \gg |Q_{n+1}|$.

 Therefore for such $\alpha$, we have $\mingap \ll N^{-1}$ for all large $N$,
 via Lemma \ref{lem 31}.
 
 \subsection{Comparison to previous results}
 For completeness, we present the constructions used in the original 
 proof of Theorem \ref{thm 61}, in terms of our notation. As mentioned above,
 the term $\frac{x+\sqrt{x^2+4\epsilon}}{2}$ is simply our $\zeta$. 
 The optional term $\sqrt{x^2+4\epsilon}$ belongs to $\sqrt{D}\cdot\Q^{\times}$.
 Therefore the numbers in $K$ covered by Theorem \ref{thm 61} are those of the
 form 
 $$ \alpha(\zeta; r,a,b) = r\cdot\zeta^a\cdot\sqrt{D}^b $$
 where $\zeta \in \OK^\times$, $r \in \Q^\times$, $a \in \Z$ and $b = 0,1$.
 The term $r = c/d$ is dealt with by multiplying all numerators by $c$ and denominators
 by $d$, and we reduce to the case $r = 1$ (and thus $\alpha \in \OK$).
 These $\alpha$ have the special property that the sequences of traces
 $\tr(\alpha\zeta^N) = \tr(\sqrt{D}^b \zeta^{N+a})$, and $\tr(\zeta^N)$ (as well
 as $\sktr(\alpha\zeta^N)$ and  $\sktr(\zeta^N)$) can have useful
 divisibility properties, for good choices of $N$, and their ratios
 provide good approximations to $\alpha$. For general $\alpha$, we can
 generalize this method and generate good approximations to $\alpha$
 as ratios of the sequences of $\tr(\alpha\omega\zeta^N)$ and $\tr(\omega\zeta^N)$,
 for any $\omega \in \OK$. The quality of the approximation is slightly worse:
 the upper bound on the error grows by the norm of $\omega$, but this is 
 not a problem if $\omega$ is fixed and $N$ is large. In order to have good
 divisibility properties for these sequences we want $\alpha\omega\zeta^N$ and
 $\omega\zeta^N$ to be $L$-th and $L'$-th powers, respectively. We achieve this
 by setting $\omega = \alpha^M$ and choosing appropriate values of $M$ and $N$ --
 but other choices for $\omega$ and $N$ might also be possible. 

\bibliography{quad_irrat_approx}{}
\bibliographystyle{hacm}

 \end{document}